\documentclass[twoside]{article}
\usepackage{graphicx,amssymb,mathrsfs,amsmath}
\catcode`@=11
\long\def\@makefntext#1{\noindent #1}
\newskip\tabcentering \tabcentering=1000pt plus 1000pt minus 1000pt
\def\REF#1{\par\hangindent\parindent\indent\llap{#1\enspace}\ignorespaces} 
\def\MCH#1#2{\setbox0=\hbox{\raise#1\hbox{#2}}\smash{\box0}}

\def\@evenfoot{}\def\@oddfoot{}

\def\@evenhead{\hbox to\textwidth{\footnotesize\rm\thepage \hfill
{\it Congdian Cheng}}} 

\def\@oddhead{\hbox to \textwidth{\footnotesize{\it
 Second Moment with A Class of
Stochastic Completion Time} \hfill\thepage}}

\def\sec#1{\vskip 3mm\leftline{\bf #1}\vskip 1mm}

\def\th#1{\vskip 1mm\noindent{\bf #1}\quad}

\def\proof{\vskip 1mm\noindent{\it Proof}\quad}

\catcode`@=12

\def\bc{\begin{center}}
\def\ec{\end{center}}
\def\no{\noindent}
\def\hang{\hangindent\parindent}
\def\textindent#1{\indent\llap{\qquad #1\ \ \enspace}\ignorespaces}
\def\ref{\par\hang\textindent}

\begin{document}
  \begin{titlepage}{ \begin{center}{ \large \bf
  Approach to Express The Second Moment with A Class of
Stochastic Completion Time \\ \vskip 0.2cm(Running head: Second
Moment with A Class of Stochastic Completion Time)}\\\vskip 0.4cm
 CHENG Cong Dian\\\vskip 0.1cm College of Mathematics and Systems Science,
Shenyang Normal University,
  Shenyang 110034, China\end{center}}


  \noindent{\bf Abstract:}This work presents
  an approach
to express
  the second moment of the completion time
with a preempt-repeat model job processed on a machine subject to
stochastic breakdowns by some distribution characters of the uptime,
downtime and processing time. \par\noindent \hrulefill\\\noindent CHENG Cong
Dian: College of Mathematics and Systems Science, Shenyang Normal
University,
  Shenyang, Liaoning Province 110034, People's Republic of
  China.\\\noindent E-mail: zhiyang918@163.com.\par\noindent
\hrulefill\\\noindent\\
\noindent {\large \bf Mathematical Subject Classification (2000)
90B36(Primary); 60K99\linebreak(Secondary)}
\end{titlepage}
\abovedisplayskip=6pt plus 1pt minus 1pt \belowdisplayskip=6pt plus
1pt minus 1pt
\thispagestyle{empty} \vspace*{-1.0truecm} \noindent
\vskip 10mm

\bc{\large\bf Approach to Express The Second Moment with A Class of
Stochastic Completion Time\\[2mm]

\footnotetext{\footnotesize Received }
} \ec
\vskip 5mm \bc{\bf CHENG Cong Dian }\\

{\small\it  College of Mathematics and Systems Science, Shenyang
Normal University,\\
  Shenyang 110034, China \\ E-mail:
  zhiyang918@163.com}\ec

\vskip 1 mm

\noindent{\small {\small\bf Abstract} This work presents
  an approach
to express
  the second moment of the completion time
with a preempt-repeat model job processed on a machine subject to
stochastic breakdowns by some distribution characters of the uptime,
downtime and processing time.\\

\vspace{1mm}\baselineskip 12pt

\no{\small\bf Keywords}   Second moment, completion time, machine,
stochastic breakdowns
   \\

\no{\small\bf MR(2000) Subject Classification} 90B36, 60K99\\
{\rm }}

 \no

\sec{1\quad Introduction}

\no Glazebrook$^{[1]}$ (1984) investigated a class of scheduling
problem with machine subject to breakdowns.
In 1990,  Birge et al.$^{[2]}$
 provided a description of the processing environment with
    machine subject to a sequence stochastic breakdowns in the sense
    as follows.

    For $k=1,2,\cdots$, the $k$-th period of machine uptime and downtime
    is characterized by two nonnegative random variables $U_k$ and $D_k$. $U_k$
    represents the $k$-th machine uptime, namely the length of the
    period between the $(k-1)$-th and $k$-th breakdown. $D_k$
    represents the $k$-th machine downtime, namely the length of the
     $k$-th breakdown. Further, $\{U_k\}$ and $\{D_k\}$ are  the
     sequences of
     independent and identically distributed (i.i.d.) nonnegative random
     variables.
     Moreover, the
     uptimes are independent of the downtimes.

      We call the process
     environment as BP Environment next.

      Note that the machine can be
     interpreted as a dynamic system, such as a person, a market or a workshop;
      and the machine uptimes and
     downtimes can be interpreted as two kinds of different states,
     such as healthful and diseased, prosperous and desolate or robust state
      and non-robust.
     The BP Environment
     can  be interpreted as a alternative system of two
     different states. Therefore, we can easily know that a lot of
     real situation can be modeled as the BP Environment.
       For example,  working on line can be interpreted as a BP Environment for it is
 inevitably disturbed
   by some unexpected events, such as viruses,  damages of component,
  power cuts, etc., randomly.
  Due to the health and mood of a person being change, we can divide
  the state of a person into two states,  healthful  state and diseased
  state, in terms of a threshold determined by theory and experience. That
  is, a person can be regarded as health (uptime) if his state satisfies the threshold,
   otherwise as disease (downtime).  According to the same way, a economic
   system
   can also be interpreted as a BP Environment. Moreover,
 the traffic of cities, which is randomly blocked up, is either
  an important instance of BP Environment.
In a word, the BP Environment is very important in practice, with
which the stochastic scheduling problems have been extensively
studied for recent twenty years, e.g. see [2-9].

For  BP Environments, a job to be processed under a BP Environment
is called  as a preempt-repeat model job if the work done on the job
is lost when a breakdown happen before it is completed, that is it
must be processed from beginning after a breakdown except it has
been completed before the breakdown happen, see [2]. The situation
to process a repeat model job under a BP Environment can be found in
many industrial application. It arise commonly in process
industries, where the product requires continuous being processed
without interruption. One example is metal refinery, where the raw
materials are purified by melting them in very high temperature. If
a breakdown (such as power outage) occurs before the metal is
purified to the required level, it will quickly cool down and the
heating process has to start again after the breakdown is over.
Other examples include running a program on a computer, downloading
a file from the internet, implementing a reliability test on a
facility, etc. Generally, if a job must be continuously processed
with no interruption until it is  completed totally, then the
processing patten of the job should be  modeled by the
preempt-repeat formulation in the presence of machine breakdowns.

Given a BP Environment, let $j$ be a repeat version job to process
under the BP Environment and with processing time $\mathbf{p}_k$ for
$k$-th  uptime, and $\{\mathbf{p}_k\}$  be an i.i.d. sequence of
random  variables. Then, it is obvious that the real processing time
$R$ of the job is a random  variable. For the random  variable $R$,
one of important topic is to express its first and second moments
$E[R]$ and $E[R^2]$. Under the conditions that $U_k$ follows an
exponential distribution, $D_k$ is a continuous random  variable
with first moment and $\mathbf{p}_k=t$ for any $k$, Birge et
al.$^{[2]}$ presented an expression of $E[R]$ with some distribution
characters of the uptime and downtime. In the present work, we
address the second moment of $R$. We will propose an approach to
express $E[R^2]$ with some distribution characters of the uptime,
downtime and processing time under the conditions that $U_k, D_k$
and $\mathbf{p}_k$ are general nonnegative continuous random
variables.

The rest of the paper is arranged as follows. Section 2 presents the
preliminary formulations. Section 3 is  devoted to express the
second moment. Section 4 provides a few of remarks. Section 5 make a
simple conclusion.

 \sec{2\quad
Formulation} This section provides a further description of the
problem  we will attack and establishes necessary formulations for
the sequel researchers.

Throughout the paper, all variables and constants are real-valued.
And variables denoted by large letters or black body letters
represent random variables. For example, $U_k$ expresses that the
$k$-th uptime is a random variable; $\mathbf{p}_j$ expresses that
the processing time of job $j$ is a random variable.

We call the problem that we will attack to be $\mathbf{SMCT}$
(second moment of the completion time), its meaning  is as follows.

 \th{SMCT.} For a given BP Environment, process  a preempt-repeat model job $j$  on
 the machine. Assume: (a) $\{U_k\}$, $\{D_k\}$ and $\{\mathbf{p}_k\}$ are
mutually independent, and all the sequences of i.i.d. and
nonnegative continuous random variables; (b) $D_k$ is
 with finite first and second
moment $\mu$ resp. $\nu$; (c) the machine is available at the
beginning of process, that is, the uptimes and downtimes are
arranged as: $U_1, D_1, U_2, D_2, \cdots, U_k, D_k, \cdots$; (d) The
starting time of the job  ( the time when $j$ begin to process ) is
zero. Let $R$ be the real completion time. Express $E[R]$ and
$E[R^2]$ with some distribution characters of the uptime, downtime
and processing time.

For convenience afterwards, we now construct the following formulae
expressions.

Let $I_A$ denote the indicator of an event A, which takes value 1 if
A occurs and otherwise 0. In terms of the meaning of $R$, if
$R<+\infty$, then we have
$$\begin{array}{rcl}
R & = & \mathbf{p}_1\cdot
I_{\{U_1\geq\mathbf{p}_1\}}+[\mathbf{p}_2+(U_1+D_1)]\cdot
I_{\{U_1<\mathbf{p}_1, U_2\geq\mathbf{p}_2\}}+\cdots\\
  &   &+[\mathbf{p}_n+\sum\limits_{k=1}^{n-1}(U_k+D_k)]\cdot
I_{\{U_i<\mathbf{p}_i, U_n\geq\mathbf{p}_n: 1\leq i\leq n-1\}}+\cdots\\
  & = & \mathbf{p}_1\cdot I_{\{U_1\geq\mathbf{p}_1\}} +
  \sum\limits_{n=2}^\infty[\mathbf{p}_n+\sum\limits_{k=1}^{n-1}(U_k+D_k)]\cdot
I_{\{U_i<\mathbf{p}_i, U_n\geq\mathbf{p}_n: 1\leq i\leq n-1\}}.
\end{array}\eqno (1)$$
Note  that
$$\begin{array}{rl}
  & I_{\{U_i<\mathbf{p}_i,U_n\geq\mathbf{p}_n;1\leq i\leq n-1\}}
  \cdot I_{\{U_i<\mathbf{p}_i,U_m\geq\mathbf{p}_m:1\leq i\leq
  m-1\}}\\
= &\left\{\begin{array}{ll}0, m\neq n\\
I_{\{U_i<\mathbf{p}_i,U_n\geq\mathbf{p}_n:1\leq i\leq n-1\}}, m=n.
\end{array}\right.
\end{array}$$
When $R<+\infty$, (1) leads to
$$\begin{array}{rcl}
R^2 & = & \mathbf{p}_1^2\cdot
I_{\{U_1\geq\mathbf{p}_1\}}+[\mathbf{p}_2+(U_1+D_1)]^2\cdot
I_{\{U_1<\mathbf{p}_1,U_2\geq\mathbf{p}_2\}}+\\
  &   &\sum\limits_{n=3}^\infty[\mathbf{p}_n+\sum\limits_{k=1}^{n-1}(U_k+D_k)]^2\cdot
I_{\{U_i<\mathbf{p}_i,U_n\geq\mathbf{p}_n:1\leq i\leq n-1\}}.
\end{array}\eqno (2)$$
\par
Finally,  set $a=E[\mathbf{p}_k|U_k\geq \mathbf{p}_k],
b=E[U_k|U_k<\mathbf{p}_k], c=E[\mathbf{p}_k^2|U_k\geq \mathbf{p}_k],
d=E[U_k^2|U_k<\mathbf{p}_k]$ and $q=P\{U_k<\mathbf{p}_k\}$ for every
$k\geq 1$.

\sec{3\quad Approach} This section focuses on expressing $E[R]$ and
$E[R^2]$.

 \th{Theorem 1.}{\it For a given
problem \textbf{SMCT}, if $0<q<1$, then
$$
E[R]=a+(b+\mu)\cdot \frac{q}{1-q}\hspace*{2mm}.\eqno (3)
$$}
\proof It is obvious that
$\{R<+\infty\}=[\{U_1\geq\mathbf{p}_1\}\bigcup(\bigcup\limits_{n=2}^\infty\{U_i<\mathbf{p}_i,
U_n\geq\mathbf{p}_n: 1\leq i\leq n-1\})]$. Note that
$\{U_i<\mathbf{p}_i, U_n\geq\mathbf{p}_n: 1\leq i\leq
n-1\}\bigcap\{U_i<\mathbf{p}_i, U_m\geq\mathbf{p}_m: 1\leq i\leq
m-1\}=\varnothing$ when $n\neq m$. We have
$$\begin{array}{rcl}P\{R<+\infty\}&=&P\{U_1\geq\mathbf{p}_1\}+\sum\limits_{n=2}^\infty P\{U_i<\mathbf{p}_i,
U_n\geq\mathbf{p}_n: 1\leq i\leq
n-1\}\\&=&\sum\limits_{n=1}^\infty q^{n-1}(1-q) =1.\end{array}$$
This shows (1) holds almost surely. Hence, we obtain
$$\begin{array}{rcl}
E[R] & = &E[\mathbf{p}_1\cdot
I_{\{U_1\geq\mathbf{p}_1\}}]+\sum\limits_{n=2}^\infty E
[\mathbf{p}_n\cdot
I_{\{U_i<\mathbf{p}_i,U_n\geq\mathbf{p}_n:1\leq i\leq n-1\}}]\\
  &  & +\sum\limits_{n=2}^\infty E[(\sum\limits_{k=1}^{n-1}D_k)
  \cdot I_{\{U_i<\mathbf{p}_i,U_n\geq\mathbf{p}_n:1\leq i\leq n-1\}}]\\
  &  &+\sum\limits_{n=2}^\infty \sum\limits_{k=1}^{n-1}E[U_k\cdot
I_{\{U_i<\mathbf{p}_i,U_n\geq\mathbf{p}_n:1\leq i\leq n-1\}}].
\end{array}\eqno (4)$$
Note that {$\{\mathbf{p}_k\}$}, {$\{U_k\}$} and {$\{D_k\}$} are
mutually independent. From (4), we get
$$\begin{array}{rcl}
E[R] & = &E[\mathbf{p}_1|U_1\geq
\mathbf{p}_1]\cdot\mathbf{p}\{U_1\geq \mathbf{p}_1\}+
\sum\limits_{n=2}^\infty E [\mathbf{p}_n\cdot
I_{\{U_n\geq\mathbf{p}_n\}}]\cdot\\
 &  &  (\prod\limits_{i=1}^{n-1}P\{U_i<\mathbf{p}_i\})
 +\sum\limits_{n=2}^\infty (n-1)\mu
  (\prod\limits_{i=1}^{n-1}P\{U_i<\mathbf{p}_i\}) \cdot
  P\{U_n\geq\mathbf{p}_n\}\\
 &  & +\sum\limits_{n=2}^\infty \sum\limits_{k=1}^{n-1}E[U_k\cdot
I_{\{U_k<\mathbf{p}_k\}}]\cdot(\prod\limits_{i=1,i\neq k
}^{n-1}P\{U_i<\mathbf{p}_i\})\cdot P\{U_n\geq \mathbf{p}_n\}\\
 & = & a(1-q)+\sum\limits_{n=2}^\infty a(1-q)q^{n-1}+
 \sum\limits_{n=2}^\infty(n-1)\mu (1-q)q^{n-1}\\
 &   &+\sum\limits_{n=2}^\infty (n-1)b(1-q)q^{n-1}\\
  & = & a(1-q)+a(1-q)(\sum\limits_{n=1}^\infty q^n)+
 (b+\mu)(1-q)q(\sum\limits_{n=1}^\infty nq^{n-1})\\
   & = & a(1-q)+a(1-q)\cdot \frac{q}{1-q}+
 (b+\mu)(1-q)q\cdot(\sum\limits_{n=1}^\infty q^n)'\\
   & = & a+(b+\mu)(1-q)q\cdot (\frac{q}{1-q})'\\
   & = & a+(b+\mu)(1-q)q\cdot \frac{1}{(1-q)^2}= a+(b+\mu)\cdot \frac{q}{1-q}\hspace*{2mm}.
\end{array}$$ This completes the proof.\\\noindent
\th{Theorem 2.}{\it For a given  problem \textbf{SMCT}, if $0<q<1$,
then
$$
E[R^2]=c+(2ab+2\mu a+2\mu b+d+\nu)\cdot \frac{q}{1-q}
+2(\mu+b)^2\cdot (\frac{q}{1-q})^2.\eqno (5)
$$}
\proof Note $P\{R<+\infty\}=1$. By (2), we obtain
$$\begin{array}{rcl}
E[R^2] & = & E[\mathbf{p}_1^2\cdot
I_{\{U_1\geq\mathbf{p}_1\}}]+E[(\mathbf{p}_2+(U_1+D_1))^2\cdot
I_{\{U_1<\mathbf{p}_1,U_2\geq\mathbf{p}_2\}}]\\
  &   & +\sum\limits_{n=3}^\infty E[\mathbf{p}_n^2\cdot
       I_{\{U_i<\mathbf{p}_i,U_n\geq q_n:1\leq i\leq n-1\}}]\\
  &   & +2(\sum\limits_{n=3}^\infty E[\mathbf{p}_n(\sum\limits_{k=1}^{n-1}
       (U_k+D_k))\cdot I_{\{U_i<\mathbf{p}_i,U_n\geq\mathbf{p}_n:1\leq i\leq
       n-1\}}])\\
  &   & +\sum\limits_{n=3}^\infty
       E[(\sum\limits_{k=1}^{n-1}\sum\limits_{l=1}^{n-1}U_k\cdot U_l)\cdot
       I_{\{U_i<\mathbf{p}_i,U_n\geq\mathbf{p}_n:1\leq i\leq
       n-1\}}]\\
  &   & +2(\sum\limits_{n=3}^\infty
       E[(\sum\limits_{k=1}^{n-1}\sum\limits_{l=1}^{n-1}U_k\cdot D_l)\cdot
       I_{\{U_i<\mathbf{p}_i,U_n\geq\mathbf{p}_n:1\leq i\leq
       n-1\}}])\\
  &   & +\sum\limits_{n=3}^\infty
       E[(\sum\limits_{k=1}^{n-1}\sum\limits_{l=1}^{n-1}D_k\cdot D_l)\cdot
       I_{\{U_i<\mathbf{p}_i,U_n\geq\mathbf{p}_n:1\leq i\leq
       n-1\}}].
\end{array}\eqno (6)$$
On the other hand, we have
$$\begin{array}{rl}
     & E[(\mathbf{p}_2+(U_1+D_1))^2\cdot
I_{\{U_1<\mathbf{p}_1,U_2\geq\mathbf{p}_2\}}]\\
  = &
  E[(\mathbf{p}_2^2+2\mathbf{p}_2(U_1+D_1)+U_1^2+2U_1D_1+D_1^2)\cdot
       I_{\{U_1<\mathbf{p}_1,U_2\geq \mathbf{p}_2\}}]\\
  = & E[(\mathbf{p}_2^2\cdot I_{\{U_2\geq \mathbf{p}_2\}})\cdot
  I_{\{U_1<\mathbf{p}_1\}}]\\
    & +2E[(U_1\cdot I_{\{U_1<\mathbf{p}_1\}})\cdot
  (\mathbf{p}_2\cdot I_{\{U_2\geq \mathbf{p}_2\}})]\\
   & +2E[(D_1\cdot I_{\{U_1<\mathbf{p}_1\}})\cdot
  (\mathbf{p}_2\cdot I_{\{U_2\geq \mathbf{p}_2\}})]\\
   & +E[(U_1^2\cdot I_{\{U_1<\mathbf{p}_1\}})\cdot
  I_{\{U_2\geq \mathbf{p}_2\}})]\\
   & +2E[D_1\cdot(U_1\cdot I_{\{U_1<\mathbf{p}_1\}})\cdot
   I_{\{U_2\geq \mathbf{p}_2\}}]\\
   & +E[D_1^2\cdot I_{\{U_1<\mathbf{p}_1,U_2\geq \mathbf{p}_2\}}]\\
= & c(1-q)q+2ab(1-q)q+2\mu a(1-q)q+\\
  & d(1-q)q+2\mu b(1-q)q+\nu(1-q)q,
\end{array}\eqno (7)$$
$$\begin{array}{rl}
     & E[(\mathbf{p}_n^2\cdot
       I_{\{U_i<\mathbf{p}_i,U_n\geq\mathbf{p}_n:1\leq i\leq n-1\}}]\\
  = & E[(\mathbf{p}_n^2\cdot
       I_{\{U_n\geq\mathbf{p}_n\}})\cdot I_{\{U_i< \mathbf{p}_i:1\leq i\leq n-1\}}]\\
  = & E[(\mathbf{p}_n^2|U_n\geq \mathbf{p}_n]\cdot P\{U_n\geq
  \mathbf{p}_n\}\cdot
  (\prod\limits_{i=1}^{n-1}P\{U_i<\mathbf{p}_i\})=c(1-q)q^{n-1};
\end{array}\eqno (8)$$
$$
\sum\limits_{n=3}^\infty
c(1-q)q^{n-1}=c(1-q)(\sum\limits_{n=3}^\infty q^{n-1})=c(1-q)\cdot
\frac{q^2}{(1-q)}=cq^2,\eqno (8)'
$$
$$\begin{array}{rl}
     & E[\mathbf{p}_n(\sum\limits_{k=1}^{n-1}(U_k+D_k))\cdot
       I_{\{U_i<\mathbf{p}_i,U_n\geq\mathbf{p}_n:1\leq i\leq n-1\}}]\\
  = & \sum\limits_{k=1}^{n-1}E[(\mathbf{p}_n\cdot
       I_{\{U_n\geq \mathbf{p}_n\}})((U_k\cdot I_{\{U_k< \mathbf{p}_k\}})
        \cdot I_{\{U_i<\mathbf{p}_i:1\leq i\leq n-1,i\neq k\}}]\\
    & +D_k\cdot I_{\{U_i<\mathbf{p}_i:1\leq i\leq n-1\}})]\\
  = & \sum\limits_{k=1}^{n-1}[a(1-q)\cdot bq\cdot q^{n-2}
      +a(1-q)\cdot \mu \cdot q^{n-1}]\\
  = & (n-1)ab(1-q)q^{n-1}+(n-1)a\mu (1-q)q^{n-1};
\end{array}\eqno (9)$$
$$\begin{array}{rl}
 & \sum\limits_{n=3}^\infty[(n-1)ab(1-q)q^{n-1}+(n-1)a(1-q)\cdot \mu\cdot
q^{n-1}]\\
= & a(b+\mu)(1-q)q(\sum\limits_{n=2}^\infty
nq^{n-1})=a(b+\mu)(1-q)q(\sum\limits_{n=2}^\infty
q^{n})'\\
= & \frac{a(b+\mu)q^2(2-q)}{1-q},
\end{array}\eqno (9)'$$
$$\begin{array}{rl}
 & E[(\sum\limits_{k=1}^{n-1}\sum\limits_{l=1}^{n-1}U_kU_l)\cdot
 I_{\{U_i<\mathbf{p}_i,U_n\geq\mathbf{p}_n:1\leq i\leq n-1\}}]\\
 = & \sum\limits_{k=1}^{n-1}E[U_k^2\cdot
 I_{\{U_i<\mathbf{p}_i,U_n\geq\mathbf{p}_n:1\leq i\leq n-1\}}]+\\
   & \sum\limits_{1\leq k,l\leq n-1,k\neq l}E[U_kU_l\cdot
 I_{\{U_i<\mathbf{p}_i,U_n\geq\mathbf{p}_n:1\leq i\leq n-1\}}]\\
  = & \sum\limits_{k=1}^{n-1}E[(U_k^2\cdot I_{\{U_k<\mathbf{p}_k\}})\cdot
 I_{\{U_i<\mathbf{p}_i,U_n\geq\mathbf{p}_n:1\leq i\leq n-1,i\neq k\}}]+\\
   & \sum\limits_{1\leq k,l\leq n-1,k\neq l}E[(U_k\cdot I_{\{U_k<\mathbf{p}_k\}})
   \cdot(U_l\cdot I_{\{U_l<\mathbf{p}_l\}})\cdot
 I_{\{U_i<\mathbf{p}_i,U_n\geq\mathbf{p}_n:1\leq i\leq n-1,i\neq k,l\}}]\\
 = & (n-1)d(1-q)q^{n-1}+(n-2)(n-1)b^2(1-q)q^{n-1};
\end{array}\eqno (10) $$
$$\begin{array}{rl}
 &
 \sum\limits_{n=3}^\infty[(n-1)d(1-q)q^{n-1}+(n-2)(n-1)b^2(1-q)q^{n-1}]\\
 = & d(1-q)q(\sum\limits_{n=2}^\infty nq^{n-1})+b^2(1-q)q^2
     (\sum\limits_{n=2}^\infty n(n-1)q^{n-2})\\
 = & d(1-q)q(\sum\limits_{n=2}^\infty q^n)'+b^2(1-q)q^2
     (\sum\limits_{n=2}^\infty q^n)''\\
 = & \frac{dq^2(2-q)}{1-q}+\frac{2b^2q^2}{(1-q)^2},
\end{array}\eqno (10)' $$
$$\begin{array}{l}
 E[(\sum\limits_{k=1}^{n-1}\sum\limits_{l=1}^{n-1}U_kD_l)\cdot
 I_{\{U_i<\mathbf{p}_i,U_n\geq\mathbf{p}_n:1\leq i\leq n-1\}}]\\
 =  \sum\limits_{k=1}^{n-1}(\sum\limits_{l=1}^{n-1}E[D_k\cdot(U_k\cdot
 I_{\{U_k<\mathbf{p}_k\}})\cdot
 I_{\{U_i<\mathbf{p}_i:1\leq i\leq n-1,i\neq k\}}
 \cdot I_{\{U_n\geq \mathbf{p}_n\}}])\\
 =   \sum\limits_{k=1}^{n-1}(n-1)bq \mu q^{n-2}(1-q)=
    (n-1)^2\mu b(1-q)q^{n-1};
\end{array}\eqno (11) $$
$$\begin{array}{l}
 \sum\limits_{n=3}^\infty(n-1)^2\mu b(1-q)q^{n-1}\\
 =\mu b(1-q)[\sum\limits_{n=3}^\infty(n-1)(n-2)q^{n-1}+
 \sum\limits_{n=3}^\infty(n-1)q^{n-1}]\\
 = \mu b(1-q)[(\sum\limits_{n=2}^\infty n(n-1)q^{n-2})q^2+
( \sum\limits_{n=2}^\infty nq^{n-1})q]\\
 = \mu b(1-q)[(\sum\limits_{n=2}^\infty q^n)''q^2+
( \sum\limits_{n=2}^\infty q^n)'q]\\
= \frac{\mu b(2-q)q^2}{1-q}+\frac{2\mu bq^2}{(1-q)^2},
\end{array}\eqno (11)' $$
$$\begin{array}{l}
 E[\sum\limits_{k=1}^{n-1}\sum\limits_{l=1}^{n-1}D_kD_l\cdot
  I_{\{U_i<\mathbf{p}_i,U_n\geq\mathbf{p}_n:1\leq i\leq n-1\}}]\\
 = (\sum\limits_{k=1}^{n-1}E[D_k^2]+
 \sum\limits_{1\leq k,l\leq n-1,k\neq l}E[D_kD_l])\cdot
 q^{n-1}(1-q)\\
= [(n-1)\nu+(n-1)(n-2)\mu^2](1-q)q^{n-1};
\end{array}\eqno (12) $$
$$\begin{array}{l}
 \sum\limits_{n=3}^\infty[(n-1)\nu+(n-1)(n-2)\mu^2](1-q)q^{n-1}\\
 =\nu(1-q)( \sum\limits_{n=3}^\infty(n-1)q^{n-1})+
 \mu^2(1-q)(\sum\limits_{n=3}^\infty(n-1)(n-2)q^{n-1})\\
 =\nu(1-q)q( \sum\limits_{n=2}^\infty q^n)'+
 \mu^2(1-q)q^2(\sum\limits_{n=2}^\infty q^n)''\\
 =\frac{\nu(2-q)q^2}{1-q}+\frac{2\mu^2q^2}{(1-q)^2}.
\end{array}\eqno (12)' $$
Combining (6)-(12), we obtain
$$\begin{array}{rcl}
 E[R^2]&=& c(1-q)+(c+2ab+2\mu a+2\mu b+d+\nu)(1-q)q+\\
       & & cq^2+\frac{2a(b+\mu)q^2(2-q)}{1-q}+[\frac{dq^2(2-q)}{1-q}
           +\frac{2b^2q^2}{1-q}]+\\
       & & 2[\frac{\mu b(2-q)q^2}{1-q}+\frac{2\mu bq^2}{(1-q)^2}]+
           [\frac{\mu (2-q)q^2}{1-q}+\frac{2\mu^2q^2}{(1-q)^2}]\\
       &=& c[(1-q)+(1-q)q+q^2]+2a[(b+\mu)(1-q)q+
           \frac{(b+\mu)(2-q)q^2}{1-q}]+\\
       & & d[(1-q)q+\frac{(2-q)q^2}{1-q}]+\nu[(1-q)q+
           \frac{(2-q)q^2}{1-q}]+\\
       & & 2\mu b[(1-q)q+\frac{(2-q)q^2}{1-q}+
           \frac{2q^2}{(1-q)^2}]+\frac{2b^2q^2}{(1-q)^2}
           +\frac{2\mu^2q^2}{(1-q)^2}\\
       &=& c+[2ab+2\mu(a+b)+d+\nu]\cdot\frac{q}{1-q}
           +2(\mu+b)^2\cdot(\frac{q}{1-q})^2.
\end{array} $$The Proof completes.
\sec{4\quad Remarks} (i) For the Theorem 1, actually $D_k$ needs
only to possess first moment.\par
 (ii) When $q=0$, we have $P\{U_1\geq \mathbf{p}_1\}=1$. Thus
 $R=\mathbf{p}_1$ almost surely. Furthermore,
 $E[R]=E[\mathbf{p}_1]$ and $E[R^2]=E[\mathbf{p}_1^2]$. When
 $q=1,P\{U_k<\mathbf{p}_k\}=1$ for all $k\geq 1$. So $j$ cannot be
 completed almost surely in the case.
 \par
 (iii) Suppose that $U_k$ follows nonnegative exponential distribution with
 parameter $\lambda$ and $ \mathbf{p}_k=t$ for all  $k\geq 1 $. Then
$$\begin{array}{l}
a=E[t|U_1\geq t]=t,\\
 b=E[U_1|U_1<t]=\frac{\int_0^t \lambda se^{-\lambda
s}ds}{P\{U_1<t\}}=\frac{1}{\lambda}-\frac{t}{e^{\lambda t}-1},
and\\
q=P\{U_1<t\}=\int_0^t \lambda e^{-\lambda s}ds=1-e^{-\lambda t}.\\
\end{array} $$
Substituting $a, b$ and $q$ with these results in  (3)
respectively, we obtain
$$
E[R]=( \frac{1}{\lambda}+\mu)(e^{\lambda t}-1).
$$
This is the same with the result of Birge et al. [2].
\par
 (iv) When $\mu =o(1)$ and $\nu =o(1)$, we call the Birge Processing Environment
as Instantaneous Birge Processing Environment (i.e., Birge
Processing Environments with stochastic instantaneous breakdowns).
Under this situation, we have
$$\begin{array}{rcl}
E[R]&=&a+\frac{bq}{1-q}+o(1),\\
E[R^2]&=&c+(2ab+d)\cdot \frac{q}{1-q}+2b^2(\frac{q}{1-q})^2+o(1).
\end{array}$$
In practice, for large part conditions, $\mu$ and $\nu$ are
actually small enough to processing time $ \mathbf{p}$ and uptime
$U_k$. This is often the realistic case. That is, the
Instantaneous Birge Processing Environment is very important in
practice.

\sec{5\quad Conclusion}
   We have successfully express the second moment of the completion time about
    a preempt-repeat model job processed on a machine subject to
    stochastic breakdowns, by some distribution characters of the uptime, the downtime
     and the processing time. The result and the method we use  are going to
     largely stimulate the development of the
    research
    on the  problems with machines
    subject to stochastic breakdowns.

{\bf Acknowledgement}
 \par
  The author
would like to thank the anonymous referees for their
 valuable suggestions and comments.
 \vskip0.2in \no {\bf References}
\vskip0.1in
 \footnotesize

 \REF{[1]\ } Glazebrook K D. Scheduling stochastic jobs on a single
machine subject to  breakdowns. \emph{Naval
 Res Logist},  {\bf31}: 251--264 (1984)

\REF{[2]} Birge J R,  Frenk J B G,  Mittenthal,
    J, et al.  Single-machine scheduling subject to stochastic breakdowns.
\emph{Naval
 Res Logist}, {\bf37}:   661-677 (1990)
\REF{[3]} Mittenthal J, Raghavachair M.   Stochastic single machine
scheduling with quadratic early-tardy penalties. \emph{Oper.Res.},
{\bf41}(4):
 786-796 (1993)

 \REF{[4]} Qi X D,  Yin  G, Birge, J R.
Single-machine scheduling with random machine breakdowns and
randomly compressible processing times.
\emph{Stochastic.Anal.Appl.}, {\bf18}(6):  635-653 (2000)

\REF{[5]\ }Jia  C F. Stochastic single machine scheduling with an
exponentially distributed due date. \emph{Operations Research
Letter} {\bf28}: 199--203  (2001)

\REF{[6]}  Jia C F.  Scheduling with random processing times to
minimize completion time variance on a single machine subject to
stochastic breakdowns. Proceedings of the 4th World Congress on
Intelligent Control and Automation:
  743-747 (2002)

\REF{[7]} Cai X, Wu X,  Zhou X. Dynamically optimal policies for
stochastic scheduling subject to preemptive-repeat machine
breakdowns. Ieee Transactions on Automation Science and Engineering,
{\bf2}(2): 158-172 (2005)

\REF{[8]\ }Tang H Y, Zhao C L, Cheng C D.   Single machine
stochastic JIT scheduling problem subject to machine breakdowns. Sci
China Ser A-Math 51:  273-292 (2008)

 \REF{[9]\ }Cheng C D, Tang H Y, Zhao C L. Scheduling jobs on a
machine  subject to stochastic breakdowns to minimize absolute
early-tardy penalties.  Sci China Ser A-Math 51:  864--888 (2008)

\end{document}